\documentclass{article}
\usepackage[T1]{fontenc}
\usepackage[latin1]{inputenc}

\makeatletter

\providecommand{\LyX}{L\kern-.1667em\lower.25em\hbox{Y}\kern-.125emX\@}

\makeatother

\begin{document}

An Extension to Fermat's Factorisation and a simple primality test\\
\\
\\
\\
Satyabrata Adhikari\\
Calcutta Mathematical Society \\
AE-374, BidhanNagar\\
Kolkata-700 064\\
\\
Abhijit Sen\\
Theory Division,Theory Group\\
Saha Institute of Nuclear Physics\\
1/AF,BidhanNagar\\
Kolkata-700 064\\
\\

Abstract\\
\\
\( \;  \)An extension to the factorisation principle as suggested by Fermat
is presented.We start from a symmetry of natural numbers and obtain the factorisation
principle therefrom.Later it is extended further to test the primality of any
natural number and finally used to factorise any given number.\\
\\
\\
\( \;  \)Factorisation of natural numbers have been a subject of interest since
the time of Euclid {[}1{]}. Subsequently many factorisation algorithms have
been devised {[}2,5{]}. In the present piece of work we try to obtain one such
factorisation algorithm as suggested by Fermat {[} 6 {]} using a symmetry of
odd composites {[} 3 {]} and then extend it to achieve further. \\
\( \;  \) Let us begin with a number N which is to be factorized.The number
N can be either even or odd.If it be even then we successively divide it by
2 to obtain an odd number.Let us suppose that after m divisions by 2 we get
an odd number P. Thus \\
\\
N = 2\( ^{m} \) {*} P .............................................................................(1)\\
\\
\( \;  \)So,we are left with splitting the odd number P into factors if possible
i.e if P be composite,else we cannot split P further into smaller factors.We
assume that P is not a perfect square.Any number that is not a perfect square
would henceforth be referred to as a non-square number.We check for it and if
so,reduce P into products of non-square terms.We now continue with each such
factor. For the sake of discussions we retain our notation P,however it is emphasized
here that this new P would in practice stand for any non-square factor of the
original P. The test if a number is a perfect square can be made simpler by
the fact that all perfect squares (expressed in base-10) have the following
22 possibilities of the last two digits {[} 4 {]}:\\
00 , e1 , e4 , 25 , o6 and e9 where e denotes even numerals 0,2,4,6,8 and o
denotes the odd numerals 1,3,5,7,9. All other combinations can be immidiately
set as being non-square .\\
\\
\( \;  \)Now it is known that {[}3{]} any odd composite number P can be expressed
in the form \\
 \\
P = ( 2n + r )\( ^{2} \) - ( r - 1 )\( ^{2} \) ...................................................(2)\\
for some natural n and r.\\
\\
Expressing (2) in the form of a quadratic equation in n we have that the discriminant
D ought to be a perfect square. Now we recast \( \sqrt{D} \) in the form \\
\\
\( \sqrt{D} \) = 4 \( (r^{2}-2r+A \))\( ^{\frac{1}{2}} \)= 4\{(r\( ^{2} \)-2r-a)
+ (A+a)\}\( ^{\frac{1}{2}} \)= (t\( _{1} \)+ t\( _{2} \))\( ^{\frac{1}{2}} \)....
( 3 )\\
 where \\
 a ) A = P + 1 \& \\
 b ) 'a' is a non-negative integer such that t\( _{1} \)= 0 and t\( _{2} \)
is a perfect square where\\
 t\( _{1} \) = (r\( ^{2} \)-2r-a ) and\\
 t\( _{2} \) = (a + A ) \\
 \\
The certainity that such an 'a' can always be found for a composite P would
be justified shortly.\\
Equating t\( _{1} \) = 0 we get for the discriminant of the quadratic equation
in r \\
\\
D = 4{*}( 1 + a ) which needs to be a perfect square since r is natural i.e\\
 ( 1 + a ) = c\( ^{2} \).......................................................................................
( 4 )\\
\\
Also from ( 3 ) we have\\
 ( A + a ) = b\( ^{2} \)......................................................................................(
5 ) \\
\\
From ( 4 ) and ( 5 ) we get \\
 b\( ^{2} \)- c\( ^{2} \) = P ...................................................................(
6 )\\
\\
It is to be noted here that (6) holds for all P,prime or composite.For a prime
P it holds only for c= b - 1. However if P be composite then there definitely
exists at least one separate pair of b and c for which c\( \neq b \)-1 as given
by ( 2 ) .This is due to the fact that ( 2 ) spans the set of all odd composites
{[} 3 {]}. This justifies the statement that a suitable 'a' can always be found.\\
Now taking successively b =ceil ( \( \sqrt{P} \) ) ,ceil (\( \sqrt{P} \) +
1) ,... etc. and for each step calculating ( b\( ^{2} \)- P ) we find a value
of b = b\( _{opt} \) for which ( b\( ^{2} \)\( _{opt} \) - P ) is a perfect
square = c\( ^{2} \)( say ). Now, to ascertain if ( b\( ^{2} \) - P ) is a
perfect square for any b we need not always find the square root explicitly
as we saw earlier.As evident from the prior discussions we do not accept values
of b,c for which c=b-1. \\
\\
If we can obtain some 'b' and 'c' with c \( \neq  \) ( b - 1 ) then we have
\\
 P = ( b + c ) {*} ( b - c ) which is Fermat's factorisation{[} 6 {]}. Otherwise
P is definitely prime and cannot be further factorised except for the trivial
result \\
P= P{*}1.\\
 Let us now try to impose some restrictions on b and c. \\
\\
We have a + A = b\( ^{2} \)\\
which clearly indicates b \( \geq  \) \( \sqrt{A} \)...................................
( 7 )\\
Now recalling that\\
\\
b{*}c = (\( \frac{b+c}{2} \))\( ^{2} \) - (\( \frac{b-c}{2} \))\( ^{2} \),
we note that the maximum value of b and c occur if c = ( b - 1 ) = ( k - 1 )
, (say) i.e \\
 k\( ^{2} \)- ( k -1 )\( ^{2} \) = ( 2 {*} k - 1 ) = P ...............................................................
(8 )\\

Therefore, b \( \leq  \) k which implies using ( 7 ) \\
 \( \sqrt{A} \) \( \leq  \)b \( \leq  \) k...........................................................................(
9 ) \\
\\
Now keeping in mind that 'c' is necessarily natural we get from ( 8 ) that \\
1 \( \leq  \)c \( \leq  \) ( k - 1 ) ............................................................(
10 )\\
\\
Had we ignored the non-square condition of P we would have required to have
\\
 0 \( \leq  \) c \( \leq  \) ( k - 1 ) ...........................................................(
10a )\\
\\
Test for primality:\\
\\
Let N ( > 2) be odd, i.e m = 0 and N = P.\\
(We safely assume this as apart from '2', all evens are composites)\\
If no set of ( b , c ) can be found satisfying ( 6 ),( 9 ) and ( 10 ) then we
can say that N is prime.\\
\\
Factorisation : \\
\\
Given a number N we do the following to factorise it:\\
\\
1. For an even N, separate out the '2' factors using representation ( 1 ).\\
( If N be odd we overlook this step in which case P = N)\\
2. Take the obtained P and reduce it to a product of non-square terms (if necessary)
and then for each such factor using the technique as above factorise it.\\
3. Take each factor and repeat the steps as done for P in ( 2 ) above.\\
4. Continue this process till each number is prime.\\
\\
Example: Let us factorise 176400.

We have 176400=2\( ^{4} \){*} 11025\\
Now 11025 has last two digits 25 which means it could be a perfect square.Now
taking the square root we get \\
 11025 = 105\( ^{2} \)

So we need to factorise 105.\\
Now 105 is not a perfect square.So we continue.

Therefore we look for b\( ^{2} \)- c\( ^{2} \) = 105.We start from b=ceil(\( \sqrt{105} \))=11
and see if ( 6 ) is satisfied.If not,we look with b = 11 + 1 = 12 and so on.For
b = 13 we have\\
13\( ^{2} \) - 105 = 169 - 105 = 64 = 8\( ^{2} \)which allows 105 = 21 {*}
5\\
Now we see that 21 is not a perfect square.Proceeding we get 21 =5\( ^{2} \)-2\( ^{2} \)=
7{*}3 while no b,c can be obtained for 5,7 or 3 meaning that they are prime
in nature.Therefore we have\\
 176400 = 2\( ^{4} \) {*} 105\( ^{2} \) = 2\( ^{4} \) {*} 3\( ^{2} \) {*}
5\( ^{2} \) {*} 7\( ^{2} \)\\
\\
 Conclusion: \\

The present technique of primality testing or factorisation is based on Fermat's
factorisation principle which we obtain from a symmetry of odd composites{[}
3 {]}.Furthermore setting restrictions on values of b and c should ease the
applicability of the above mentioned algorithm.

In the present article we have shown that for the factorisation of any given
number,we do not require to test the primality of it beforehand. Instead we
use the ranges of 'b' and 'c'. If we fail to find some suitable values of 'b'
and 'c' then the concerned P has to be necessarily prime. \\
\\

References:\\
\\
1.http://www.arXiv.org , math.NT/0105219\\
2.CRC Concise Encyclopedia of Mathematics,Eric W Weisstein\\
Chapman \& Hall/CRC ,Page: 1429\\
3.http://www.arXiv.org , math.NT/0109114\\
4.CRC Concise Encyclopedia of Mathematics,Eric W Weisstein\\
Chapman \& Hall/CRC ,Page: 1708\\
5. Fundamentals of Number Theory,William J. LeVeque \\
Section 5.5, Page 113, Dover 1977\\
6.CRC Concise Encyclopedia of Mathematics,Eric W Weisstein\\
Chapman \& Hall/CRC ,Page: 617
\end{document}